\documentclass[12pt, DIV12]{article}
\usepackage[utf8]{inputenc}
\usepackage[T1]{fontenc}
\usepackage{sectsty}
\sectionfont{\mdseries \centering}

\usepackage{enumerate}
\usepackage{lmodern, microtype}
 
\usepackage{titling}		 
\usepackage{authblk} 	 

\usepackage{amsmath, amsthm, amssymb}
\usepackage{aligned-overset}   

\usepackage[english]{babel}
\usepackage[nottoc]{tocbibind}

\usepackage{placeins}  
\usepackage{caption}
\usepackage{float} 

\usepackage{fancyhdr}  
\usepackage{tikz}     

\usepackage{graphicx}		
\usepackage{xcolor}				

\definecolor{darkgreen}{rgb}{0,0.5,0}
\definecolor{darkred}{rgb}{0.7,0,0}
\definecolor{darkblue}{rgb}{0,.2,.7}

\usepackage[colorlinks, 
citecolor=darkred, linkcolor=darkgreen, urlcolor=darkblue,
pdfpagelabels=true,
unicode=true
]{hyperref}

\hypersetup{
	pdfauthor={Jonah Reuß},
	pdftitle={A note on the existence of nontrivial zero modes on Riemannian manifolds},
	pdflang={en} 
}

\newtheoremstyle{definitionstyle}
{5pt}
{5pt}
{}
{}
{\bfseries}
{.}
{0.5em}
{}

\newtheoremstyle{italicstyle}
{}{} 
{\itshape} 
{} 
{\bfseries} 
{.} 
{.5em} 
{} 

\theoremstyle{italicstyle}
\newtheorem{satz}{Satz}
\newtheorem{thm}[satz]{Theorem}
\newtheorem{defi}[satz]{Definition}
\newtheorem{prop}[satz]{Proposition}

\theoremstyle{definitionstyle}

\newtheorem{rem}[satz]{Remark}
\newtheorem{ex}[satz]{Example}

\newcommand{\psieps}{|\psi|_{\varepsilon}}
\newcommand{\cl}{\cdot_{\text{Cl}}} 
\newcommand{\repsone}{R^{(1)}_{\varepsilon}} 
\newcommand{\repstwo}{R^{(2)}_{\varepsilon}}
\newcommand{\peps}{P_{\varepsilon}}

\title{\textbf{A note on the existence of nontrivial zero modes on Riemannian manifolds}}
\author{Jonah Reu\ss}
\affil{Mathematisches Institut, Universität Freiburg, 79100 Freiburg, Germany  \\
E-mail: \texttt{jonah.reuss@math.uni-freiburg.de}}

\begin{document}
	\setcounter{page}{0}
	\pagestyle{fancy}  
	\pagenumbering{arabic}  
	\date{}
	\maketitle
	\begin{abstract}
		\noindent We prove a necessary criterion for the (non-)existence of nontrivial solutions to the Dirac equation $D\psi= i A\cl \psi$ on Riemannian manifolds that are either closed or of bounded geometry. This generalizes a result of Rupert Frank and Michael Loss on $\mathbb{R}^n$ where the criterion relates the $L^n$-norm of $A$ to the Sobolev constant on $\mathbb{R}^n$. On Riemannian manifolds the role of the Sobolev constant will be replaced by the Yamabe invariant. If $n$ is odd, we show that our criterion is sharp on $\mathbb{S}^n$. 
	\end{abstract}

	\section{\textsc{Introduction}}

	\renewcommand\thefootnote{}
	\footnotetext{2020 \textit{Mathematical subject classification}. Primary 53C27, Secondary 34L40, 53C18\\
	Keywords: Dirac operator, spin geometry, conformal geometry}
	\renewcommand\thefootnote{\arabic{footnote}}
	In \cite[Theorem 1]{FrankLoss2} R. Frank and M. Loss prove a sharp necessary condition for the (non-)existence of nontrivial solutions to the Dirac equation
	\begin{align}
		\gamma  (-i \nabla - i  A)\psi =0 \text{ in } \mathbb{R}^n \text{ for } n\geq 3.  \label{eq:Rnmaineq}
	\end{align}
	Solutions to this equation are called zero modes and were first derived in $3$ dimensions in \cite{lossyau}. They arise in several areas of mathematics and are particularly related to the size of magnetic fields in physics. We refer the reader to \cite[Section 1]{FrankLoss} and the references therein for more background information and material. In the setting in \cite{FrankLoss2} the quantity $A$ is a vector field on $\mathbb{R}^n$ and $\psi$ is a  spinor, which in this context is a $\mathbb{C}^{2^{[n/2]}}$ valued function. In the context of \eqref{eq:Rnmaineq}, the term $\gamma i A\psi$ stands for the Clifford action of $A$ on $\psi$ and $\gamma (-i\nabla)\psi$ stands for the Dirac operator on $\mathbb{R}^n$.   \\
	In this paper, we will attack the question of whether a similar statement holds on Riemannian manifolds and its associated  spinor bundles. For a vector field $X$ on $M$, we denote by $X\cl\psi$ the Clifford action of $X$ on a spinor $\psi$. Note that for $p\in M$, in contrast to the convention in \cite{FrankLoss2}, in local coordinates the Clifford algebra is determined by the multiplication rule
	\begin{align}
		e_i e_j+ e_j e_i=-2\delta_{ i j}
	\end{align}
	for an orthonormal basis $(e_i)_{i=1,...,n}$ of $T_pM$. 
	By translating \eqref{eq:Rnmaineq} into the setting of Riemannian manifolds, we investigate the existence of nontrivial solutions to 
	\begin{align}
		i A \cl \psi = D\psi \label{eq:maineq} \text{ on } M.
	\end{align} 
	Throughout, we assume that $(M,g)$ is a complete, smooth Riemannian  manifold without boundary and $\dim M=n \geq3$ together with a Clifford bundle $\Sigma$, cf. Section \ref{sec:prelimgeom} for precise definitions. Sections of the Clifford bundle are called spinors and will be mostly denoted by $\psi$. The associated Dirac operator will be denoted by $D$ as long as the metric $g$ is fixed.
	Note, that whenever we assume $M$ to be spin we assume the spin structure $\sigma$ to be fixed and $\Sigma$ will denote the associated spinor bundle to the $\text{Spin}(n)-$principal bundle over $(M,g,\sigma)$. 
	Additionally we always assume
	\begin{align*}
		A\in L^{n}(M,TM).
	\end{align*}
	The proof of \cite[Theorem 1]{FrankLoss2} relies on the Euclidean Sobolev inequality, cf. \cite[Theorem 4.4]{hebey}, in $\mathbb{R}^n$ given by
	\begin{align*}
		S_n\big(\int_{\mathbb{R}^n} |u|^{\frac{2n}{n-2}}dx\big)^{\frac{n-2}{n}} \leq \int_{\mathbb{R}^n} |\nabla u |^2dx \quad\text{for } u \in H^{1,2}(\mathbb{R}^n).
	\end{align*}
	Here, $S_n$ denotes the optimal constant in the Sobolev inequality, which was first independently derived in \cite{tal} and \cite{aub2}. Unfortunately this inequality is no longer true on general manifolds. 
	However, by \cite[p. 196]{stru} there is a close relation of
 	 the optimal constant in the Sobolev inequality on $\mathbb{R}^n$ and the so called Yamabe invariant on $\mathbb{S}^n$. This suggests to employ an inequality, based on the Yamabe invariant, cf. Section \ref{sec:prelim}. The Yamabe invariant, denoted by $Y(M,[g])$, is a conformal invariant, i.e., for a Riemannian metric $g$ it is independent on the choice of a metric belonging to the conformal class of $g$.\\
	It plays an important role in the solution of the original Yamabe problem, which affirmatively answers the question whether there always exists a metric of constant scalar curvature that is conformally equivalent to a given metric on a closed Riemannian manifold $(M,g)$. For more details on the original Yamabe problem we refer to \cite{aub,lepa,schoe,tru,yam} and to \cite{gro} and \cite{jin} for some results on non-compact manifolds.\\
	Before stating our main results let us finally mention that the Bär-Hijazi-Lott invariant relates problems involving the Yamabe invariant and Dirac operators through the Hijazi inequality, cf. \cite{hij}. Even if this may not be directly related to the results in the present paper, it may be interesting to keep this in mind for future work to see whether it is possible to connect our results somehow with spinorial Yamabe type problems. 
	\paragraph{Statement of results}
	Our first main result provides a necessary condition for nontrivial solutions to \eqref{eq:maineq} on closed manifolds under similar conditions as in \cite[Theorem 1]{FrankLoss2}. 
	\thm[]{Let $(M,g,\sigma)$ be a closed Riemannian spin manifold, $\dim M=n\geq3$, with spinor bundle $\Sigma$. Suppose $\psi \in L^{p}(M,\Sigma)$ for some $p$ with $\frac{2n}{n-2}\leq p <\infty$ is a nontrivial solution to \eqref{eq:maineq}
		then 
		\begin{align}
			\frac{n}{4(n-1)}Y(M,[g]) \leq \|  A\|_{L^{n}}^2. \label{eq:inequYamspin}
		\end{align} \label{thm:spin}}
\rem[]{We emphasize that \eqref{eq:inequYamspin} is sharp in any odd dimension. By "sharp", we mean that in any odd dimension $n$, we can recover equality in \eqref{eq:inequYamspin} from the case in $\mathbb{R}^n$ on the standard sphere $\mathbb{S}^n$, cf. Remark \ref{rem:recovery_on_sphere}. Indeed, in Section \ref{sec:recovery} we will characterize closed manifolds on which equality in \eqref{eq:inequYamspin} is possibly attained. In three dimensions we will see that $(M,g)$ is already conformally equivalent to $(\mathbb{S}^3,g_{g_{\mathbb{S}^3}})$.}

		\rem[]{\label{rem:manandrncoincide}
		In \cite[Theorem 1]{FrankLoss2} an inequality similar to \eqref{eq:inequYamspin} is proven for $(\mathbb{R}^n,g_{\mathbb{R}^n})$, which involves the optimal constant in the Sobolev inequality on $\mathbb{R}^n$ instead of the Yamabe invariant.
		However, due to the conformal invariance of \eqref{eq:maineq} on $\mathbb{R}^n$, cf. \cite[Remark 3]{FrankLoss2}, the appearance of the Yamabe invariant in \eqref{eq:inequYamspin} is not surprising and one expects the constants on $\mathbb{R}^n$ and $\mathbb{S}^n$ to coincide, due to the conformal equivalence of $\mathbb{R}^n$ and $\mathbb{S}^n$. Indeed,
		by \cite{tal,aub2} the optimal constant in the Sobolev inequality on $\mathbb{R}^n$ is given by
		\begin{align*}
			S_n= \frac{n(n-2)}{4} |\mathbb{S}^n|^{\frac{2}{n}},
		\end{align*}
		where $|\mathbb{S}^n|$ denotes the volume of the $n$-dimensional unit sphere.
		The Yamabe invariant of the sphere is given by
		\begin{align*}
			Y(\mathbb{S}^n,[g_{\mathbb{S}^n}])= n(n-1)|\mathbb{S}^n|^{\frac{2}{n}},
		\end{align*}
		cf. \cite[p. 196]{stru}. 
		Comparing the constants in \cite[Theorem 1]{FrankLoss2} and  Theorem \ref{thm:spin} yields
		\begin{align*}
			\frac{n}{4(n-1)}Y(\mathbb{S}^n,[g_{\mathbb{S}^n}])=\frac{n^2}{4}|\mathbb{S}^n|^{\frac{2}{n}}=\frac{n}{n-2}S_n.
		\end{align*}}\\
In our second main result, we prove a similar statement for a class of non-compact manifolds, namely manifolds of bounded geometry. As some regularity issues appear we impose stronger initial regularity conditions on the solution compared to Theorem \ref{thm:spin}.
		\thm[]{Let $(M,g,\sigma)$ be a Riemannian spin manifold of bounded geometry with spinor bundle $\Sigma$, $\dim M=n\geq3$. Suppose $A\in L^{n}(M,TM)$, $\|A\|_{L^{\infty}}<\infty$ and $\psi \in L^{2}(M,\Sigma)$ is a nontrivial solution of \eqref{eq:maineq}. Then
			\begin{align}
				\frac{n}{4(n-1)}Y(M,[g]) \leq  \| A \|^2_{L^n} . \label{eq:inequYamspinbd}
			\end{align} \label{thm:mainboundedgeomspin}}
		\rem[]{It is remarkable that no additional curvature term appears in \eqref{eq:inequYamspin} and \eqref{eq:inequYamspinbd}. This is due, on the one hand, to the special geometric structure of the problem, and on the other hand, to the use of the Yamabe invariant. 
		The proof of Theorem \ref{thm:spin} relies on a certain integral identity, cf. Proposition \ref{prop:InteIdent}. In contrast to $(\mathbb{R}^n,g_{\mathbb{R}^n})$, in our integral identity, the Clifford contraction $\mathcal{K}$ of the curvature of $\Sigma$ appears. Rearranging this inequality so that the Yamabe-invariant appears leads to an extra scalar curvature term. Then inserting the explicit form of $\mathcal{K}$ in the case where $M$ is spin, cf. \cite[Proposition 4.21]{roe}, we get rid of the curvature term. This will be examined in detail in Section \ref{sec:mainthm}.
			 }
		\paragraph{Structure of the paper} We introduce notations and collect geometric and analytic preliminaries in Section \ref{sec:prelim}. In Section \ref{sec:intident}, we state and prove an integral equality, which serves as the main starting point for the proof of Theorem \ref{thm:spin}. Section \ref{sec:mainthm} is devoted to proving of a slightly more general version of Theorem \ref{thm:spin}. Here we also establish a similar theorem for manifolds of bounded geometry. In the last section we recover the case of equality on the sphere and characterize manifolds on which equality is possibly attained in \eqref{eq:inequYamspin}.
		\paragraph{Acknowledgment} I would like to thank my supervisor Nadine Gro{\ss}e for many important and enlightening discussions and suggestions about Dirac operators and spinorial problems in general. I also want to thank my second supervisor Guofang Wang for further suggestions and for bringing attention to the problem. Finally, I want to thank Rupert Frank for a very useful discussion during his visit at the University of Freiburg. 
	\section{\textsc{Preliminaries}\label{sec:prelim}}
	\subsection{Geometric preliminaries\label{sec:prelimgeom}} Let $(M,g)$ be a smooth and complete Riemannian manifold with $\dim M=n$. Let $(\Sigma \rightarrow M,\langle \cdot, \cdot \rangle, \nabla^{\Sigma})$ be a bundle of Clifford modules over $M$ equipped with a hermitian metric $\langle \cdot,\cdot \rangle$ and a metric connection $\nabla^\Sigma$. If the Clifford multiplication is skew-adjoint with respect to $\langle\cdot,\cdot \rangle$ and the connection $\nabla^{\Sigma}$ is compatible with the Levi-Civita connection on $M$ and the Clifford multiplication, we call $\Sigma$ a Clifford bundle over $M$, cf. \cite[Definition 3.4]{roe}.
	Sections of $\Sigma$ will be called \textit{spinors}. We define the Dirac operator $D_g$ associated with $\Sigma$ by the composition of the connection $\nabla^{\Sigma}$ and the Clifford multiplication, cf. \cite[Definition 3.5]{roe}. In terms of a local orthonormal frame $D_g$ is given by
	\begin{align}
		D_g \psi = \sum^n_{i=1}e_i \cl \nabla_{e_i}^{\Sigma} \psi, \label{eq:dirlocal}
	\end{align}
	where $\psi: M \rightarrow \Sigma$. If we are working with a fixed metric $g$, we will write $D=D_g$.
	 An important property of the Dirac operator is its close relation to the curvature of $\nabla^{\Sigma}$. Following \cite[Chapter 3, (3.8)]{roe} there is an operator $\mathcal{K} \in \text{End}(\Sigma)$ satisfying the Weitzenböck formula $D^2= (\nabla^{ \Sigma } )^{*}\nabla^{\Sigma} + \mathcal{K} $, where $\mathcal{K}$ is locally given by 
	\begin{align}
		\sum^n_{i<j} \left(e_j e_i\right) \cl\left(\nabla_{e_j}\nabla_{e_i}- \nabla_{e_i}\nabla_{e_j}\right) \label{eq:cliffconloc}
	\end{align}
	for a synchronous orthonormal framing $(e_i)_{i=1,...,n}$, cf \cite[Definition 2.15]{roe}. Note that at the origin of a local orthonormal synchronous frame $\nabla_{e_j}e_i=0$ and the Lie bracket of $e_i$ and $e_j$ vanishes for all $i,j$.  
	The operator $\mathcal{K}$ is called the \textit{Clifford contraction}. By \cite[Proposition 3.18]{roe}, we have the explicit formula $\mathcal{K}=F^{^\Sigma} + \frac{1}{4}\text{scal}_g$, where $\text{scal}_g$ denotes the scalar curvature of $(M,g)$ and the term $F^{\Sigma}$ is the Clifford contraction of the curvature of $\nabla^{\Sigma}$ and is known as the twisting curvature of $\Sigma$. \\
	If $M$ admits a spin structure $\sigma$, then $\Sigma(M,g,\sigma)$ is called the spinor bundle over $M$ associated to the spin structure $\sigma$. Note that whenever we talk about spinor bundles we always assume that $\sigma$ is fixed. If $g$ is fixed as well, we will simply write $\Sigma(M,g,\sigma)=\Sigma$. 
	\rem[{\cite[Proposition 4.21]{roe}}]{For a Riemannian spin manifold $(M,g,\sigma)$ with fixed spin structure $\sigma$ the twisting curvature of the spinor bundle associated to the spin structure is zero and we have $\mathcal{K}=\frac{1}{4}\text{scal}_g$. \label{rem:spin} }\\
	We will now introduce the \textit{twistor operator} $T_g$ as it will be of great importance for the proof of the main theorem and the characterization of cases of equality in \eqref{eq:inequYamspin}.
	\defi[{\cite[Section 1]{friedrich_baum}}]{The twistor operator $T_g$ is defined by the local formula
	\begin{align}
		T_g \psi = \sum_{i=1}^{n} e_i \otimes \left(\nabla_{e_i} \psi + \frac{1}{n}e_i \cl D\psi\right) \label{eq:twistor}
		\end{align}
	If $\psi \in \text{ker}(T_g)$, we call $\psi$ a twistor spinor. }
	\rem[]{According to the first remark in \cite[Section 1]{friedrich_baum}, the following are equivalent
	\begin{enumerate}
		\item[(i)] $T_g \psi =0$,
		\item[(ii)] $\nabla_X \psi + \frac{1}{n}X \cl D\psi=0$ for all  $X\in \mathfrak{X}(M)$,
		\item[(iii)] $X\cl \nabla_Y \psi + Y\cl \nabla_X \psi=\frac{2}{n}g(X,Y)D\psi$ for all $X \in \mathfrak{X}(M)$.
		\end{enumerate}
	By a computation based on the properties of the Clifford multiplication we obtain from (ii) and \eqref{eq:twistor}
	\begin{align}
		|T_g\psi|=\sum_{j=1}^{n} \Big|\big[\nabla_{e_j} + \frac{1}{n} e_j \cl D  \big]\psi\Big|^2= |\nabla \psi|^2-\frac{1}{n}|D\psi|^2,\label{eq:twistor_abs}
		\end{align}
	which will be the starting point in the proof of the main theorem in Section \ref{sec:mainthm}.\label{rem:twistor}}
	\\Lastly we introduce the notion of bounded geometry, which will be the assumption for the main theorem in the setting of non-compact manifolds.
	\defi[{\cite[Definition A.1.1]{shu}}]{Let $(M,g)$ be a Riemannian manifold of dimension $n$.
		$M$ is called a manifold of bounded geometry if the following two conditions are satisfied
		\begin{enumerate}
			\item[(i)] $r_{\text{inj}}>0$
			\item[(ii)] $|\nabla^k R|\leq C_{k},\quad k=0,1,2,...$.
		\end{enumerate}
		Here $r_{\text{inj}}$ denotes for the injectivity radius of $M$ and $\nabla^{k}R$ represents the $k$th-covariant derivative of the Riemannian curvature tensor $R$. The pointwise norm $|\nabla^{k}R |$ is norm induced by $g$. \\
		A hermitian vector bundle $S$ over $M$ is called bundle of bounded geometry if it admits a trivialization $(U_{\alpha},\kappa_{\alpha},h_\alpha)$, where $(U_{\alpha},\kappa_{\alpha})$ is a geodesic atlas, such that the transition functions $\mu_{\alpha \beta}$ are uniformly $C^{\infty}$-bounded, i.e. for any $\alpha, \ \beta$ we have on $U_{\alpha}\cap U_{\beta}$
		\begin{align*}
			|\nabla^{k}\mu_{\alpha \beta}| \leq C_{k},
		\end{align*}
		where $C_k$ is independent of $\alpha$ and $\beta$.
		 \label{def:bdgeom}}
		 \upshape
	\subsection{Analytic preliminaries}
	Let $(M,g)$ be a smooth and complete Riemannian manifold, either compact or of bounded geometry. 
	By $\Gamma_{c}^{\infty}(M,\Sigma)$, we denote the smooth and compactly supported spinors of $\Sigma$. 
	For $\psi \in \Gamma^{\infty}_c(M,\Sigma)$ we define, for $1\leq p<\infty$ the $L^{p}$-norm of $\psi$ by
	\begin{align*}
		\|\psi\|_{L^p(M,\Sigma)}:=\left(\int_M |\psi|^p dx\right)^{\frac{1}{p}},
	\end{align*}
	where $|\cdot|$ is understood as the fiber-wise norm induced by the hermitian metric and $dx$ denotes the Riemannian volume element. 
	Then for $k \in \mathbb{N}_{0}$ we define
	\begin{align*}
		\| \psi \|^p_{H^{k,p}} := \sum^k_{j=1} \|\nabla^{j} \psi \|^p_{L^p(M,\Sigma \otimes T^{*}M^{\otimes j })},
	\end{align*}
	 as the Sobolev norm of $\psi$.
	 We define the Sobolev space $H^{k,p}(M,\Sigma)$ as the closure of $ \Gamma^{\infty}_{c}(M,\Sigma)$ with respect to $\|\cdot \|_{H^{k,p}}$. Note, that there are many equivalent ways of defining Sobolev spaces on manifolds and vector bundles. Without any attempt of completeness
	  \cite[Appendix 1, (1.3)]{shu}, \cite[Section 3.1]{gro1}, \cite[Section 2.2, Section 3.1]{hebey},\cite[Chapter 5]{roe} provide good references for more details on the definition of Sobolev spaces on manifolds or vector bundles. 
	  \prop[{\cite[5.14]{roe}}]{The Dirac operator satisfies the Garding inequality on closed manifolds, i.e., there is a constant $C$ such that  for all sections $\psi$ of $\Sigma$
	  \begin{align}
	  	\|\psi \|_{H^{1,2}(M,\Sigma)} \leq C\left( \|\psi\|_{L^{2}(M,\Sigma)} + \|D\psi\|_{L^{2}(M,\Sigma)}\right). \label{eq:gard}
	  	\end{align}\label{rem:gard} }\\ 
  	\upshape
  	As the proof of \eqref{eq:gard} is reduced to a local coordinate patch by the means a partition of unity, the result still holds on manifolds of bounded geometry, cf. \cite[Section 3.1.3]{amm2}.\\
  	\\
	 For a Riemannian manifold the \textit{Yamabe invariant} is defined as 
	 \begin{align*}
	 	Y(M,[g])= \inf_{\underset{v\neq0}{v \in C_c^\infty(M)}} \frac{\int_{M} vL_g v dvol_g}{\|v\|^2_{L^{\frac{2n}{n-2}}}}, 
	 \end{align*}
	  where the operator $L_g= \frac{4(n-1)}{n-2} \Delta_g + \text{scal}_g$ is known as the conformal Laplacian. It is a conformal invariant, cf. \cite[Proposition 5.8]{aubbook}, and is one of the main objects studied in Yamabe-type problems. Following \cite[p.1]{gro}, we remark that since we take the infimum over compactly supported functions, we can use the same definition of the Yamabe invariant for non-compact manifolds.
	\section{\textsc{Integral Identity}\label{sec:intident}}
	 Analogously to the proof of \cite[Theorem 1]{FrankLoss2}, the proof of Theorem \ref{thm:spin} is based on an integrated version of the Schrödinger-Lichnerowicz formula, which will be established in this section. 
	 Throughout this section, let $(M,g)$ be a Riemannian manifold and $\Sigma$ a Clifford bundle over $M$, $\dim M= n\geq3$. Similarly as in \cite[Section 2]{FrankLoss2}, to avoid problems concerning regularity and dividing by zero we consider the following regularization for $\psi \in \Gamma(M,\Sigma)$:
	 \begin{align*}
	 	\psieps:= \sqrt{|\psi|^2+\varepsilon^2}, \qquad \varepsilon>0.
	 \end{align*}
	 We will now state and prove the integral identity analogue to \cite[Proposition 7]{FrankLoss2} and highlight the main differences to the flat case, where $\psi: \mathbb{R}^n\rightarrow\mathbb{C}^N$. By $\mathcal{K}$ we denote the Clifford contraction of the curvature of $\Sigma$. We always assume that $\mathcal{K}$ is a bounded section in the endomorphism bundle of $\Sigma$. 
	 \prop[]{If $\psi \in H^{1}(M,\Sigma)$, then for all $\varepsilon>0$
	 	\begin{align}
	 		\begin{split}
	 			&\int_{M} 
	 			\left|T_g \big(\frac{\psi}{\psieps^{n/n-1}}\big)\right|^2 \psieps^2 dx \\
	 			=\frac{n-1}{n} &\int_{M} \frac{\left|D\psi\right|^2}{\psieps^{\frac{2}{n-1}}}dx- \frac{n-1}{(n-2)^2} \int_{M} \left| \nabla^M\left(\psieps^{\frac{n-2}{n-1}}\right)\right|^2 \left[2(n-1)-n\frac{\left|\psi\right|^2}{\psieps^2}\right]dx \\ -&\int_{M}  \frac{\text{Re} \left( \left\langle \psi, \mathcal{K}\psi\right\rangle \right)}{ \psieps^{\frac{2}{n-1}} }dx. \label{eq:integraliden}
	 		\end{split}
	 	\end{align} \label{prop:InteIdent}}
 	\upshape
 	The proof follows the idea and notions of \cite[Proposition 7]{FrankLoss2}.
 	\begin{proof} As in the proof of \cite[Proposition 7]{FrankLoss2}, the proof splits into 5 steps. We first note that by standard computations we have
 	\begin{align}
 		\text{Re}\left(\left\langle \psi, \nabla^{\Sigma}\psi \right\rangle \right) =\left| \psi \right| \nabla^M \left|\psi\right| = \psieps \nabla^M \psieps.  \label{eq:step0}
 	\end{align}
 	Note that $\langle \psi, \nabla^{\Sigma}\psi \rangle$ is not the hermitian product of two spinors as $\nabla^{\Sigma} \psi \in \Gamma(M,\Sigma\otimes T^* M)$. It should be considered as a $1$-form with the local expression $\langle \psi, \nabla^{\Sigma} \psi \rangle= \sum^n_{j=1} \langle \psi, \nabla^{\Sigma}_{e_j} \psi \rangle e^b_j$. Here $e_j^b$ denotes the musical isomorphism of $T_pM$ and $T^{*}_pM$ given by $X^b(Y):=g(X,Y)$ for $X,Y\in T_pM$.\\
 	We use \eqref{eq:twistor_abs} for the integrand on the left hand side and observe
 	\begin{align*}
 			\left|T_g \big(\frac{\psi}{\psieps^{n/n-1}}\big)\right|^2 \psieps^2= \left(\left| \nabla^{\Sigma} \left(\frac{\psi}{\psieps^{n/n-1}}\right)\right|^2-\frac{1}{n}\left|D\left(\frac{\psi}{\psieps^{n/n-1}}\right) \right|^2\right) \psieps^2.
 	\end{align*}
 	Then, using \eqref{eq:step0}, step (1) and (2) in \cite[p.9, Equation (12) and (14)]{FrankLoss2} are just pointwise identities and the proof for sections in $\Gamma(M,\Sigma)$ can be adapted straightforward from the proof for sections in $\Gamma(\mathbb{R}^n,\mathbb{C}^N)$. We obtain the identities
 				\begin{align}
 				\left|\nabla^{\Sigma} \left( \frac{\psi}{\psieps^{\frac{n}{n-1}}}\right)\right|^2 \psieps^2 = \frac{\left|\nabla^{\Sigma} \psi\right|^2}{\psieps^{\frac{2}{n-1}}} + \left| \nabla^M \left( \psieps^{\frac{n-2}{n-1}} \right) \right|^2\left[ \left(\frac{n}{n-2}\right)^2\frac{\left|\psi\right|^2}{\psieps^2} - \frac{2n(n-1)}{(n-2)^2} \right] \label{eq:step1}
 			\end{align}
 			and
 			\begin{align}
 				\begin{split}
 					\left|D\left( \frac{\psi}{\psieps^{\frac{n}{n-1}}}\right)\right|^2\psieps^2=&\frac{\left|D\psi\right|^2}{\psieps^{\frac{2}{n-1}}} + \left(\frac{n}{n-2}\right)^2\left|\nabla^M \left(\psieps^{\frac{n-2}{n-1}}\right)\right|^2\frac{\left|\psi\right|^2}{\psieps^2} \\ -&\frac{2n}{n-1}\frac{1}{\psieps^{\frac{2}{n-1}+1}}\text{Re}\left( \left\langle D\psi, \nabla^M \psieps \cl \psi \right\rangle \right).
 				\end{split} \label{eq:step2}
 			\end{align}
 	For step (3) let $\chi \in C^{\infty}_{c}(M)$ such that $\chi$ is supported in a small enough coordinate patch $U\subset M$, which admits a synchronous orthonormal frame $(e_i)_{i=1,...,n}$. We claim
 	\begin{align}
 		\begin{split}
 			\int_{M} \psieps^{-\frac{2}{n-1}} \left|\nabla^{\Sigma} \psi\right|^2 \chi dx &= \frac{2(n-1)}{(n-2)^2} \int_{M} \left| \nabla^M \left( \psieps^{\frac{n-2}{n-1}}\right) \right|^2 \chi  dx \\
 			+& \int_{M} \psieps^{-\frac{2}{n-1}}\left|D \psi\right|^2 \chi dx \\
 			-& \frac{2}{n-1} \int_{M} \psieps^{-1-\frac{2}{n-1}} \text{Re} \left(\left\langle \left(\nabla^M \psieps\right) \cl \psi, D\psi \right\rangle\right) \chi dx \\
 			+&\sum_{\underset{j\neq k}{j,k=1}}^{n} \psieps^{-\frac{2}{n-1}}\text{Re}\left(\left\langle e_j \cl \psi, e_k \cl \nabla^{\Sigma}_{e_k} \psi\right\rangle\right) \nabla^M_{e_j} \chi dx \\
 			-& \int_{M} \psieps^{-\frac{2}{n-1}}\text{Re}\left(\left\langle \psi, \mathcal{K} \psi \right\rangle\right) \chi dx.
 		\end{split} \label{eq:step3}
 	\end{align}
		The proof of \eqref{eq:step3} is based on the identity
\begin{align}
	\begin{split}
		\sum^n_{k< j}&\int_{M} f \left\langle \nabla^\Sigma_{e_k} \psi, e_k e_j\cl \nabla^\Sigma_{e_j} \psi \right\rangle dx+\int_{M} f \left\langle \nabla^\Sigma_{e_j} \psi, e_j e_k\cl \nabla^\Sigma_{e_k} \psi \right\rangle dx\\=& -\sum^n_{k< j }\int_M\left( \left(\nabla^M_{e_k}f\right) \left\langle \psi, e_k e_j \cl \nabla^\Sigma_{e_j} \psi \right\rangle + \left(\nabla^M_{e_j}f\right) \left\langle \psi, e_j e_k \cl \nabla^\Sigma_{e_k} \psi \right\rangle \right)dx\\
		&- \int_{M} f \left\langle\psi,\mathcal{K}\psi\right\rangle dx.
	\end{split} \label{eq:step3integralident}
\end{align}
 	for $f=\psieps^{-2/(n-1)}\chi$. Before we prove \eqref{eq:step3integralident}, we make sure that every integral above is well defined.\\
 	As $f$ and is bounded and $\psi \in L^2(M,\Sigma)$, $\nabla\psi \in L^2(M,\Sigma \otimes T^*M)$ the left hand side and the last integral on the right hand side are well defined. For the first integral on the right hand side we note
 	\begin{align}
 		\nabla^M_{e_j} f = -\frac{2}{n-1} \psieps^{-\frac{2}{n-1}-1} \left(\nabla^M_{e_j} \psieps\right)\chi + \psieps^{-\frac{2}{n-1}}\left(\nabla^M_{e_j}\chi\right). \label{eq:step3help3}
 	\end{align}
 	Again, as $\psi \in L^2(M,\Sigma)$, $\nabla \psi \in L^2(M,\Sigma\otimes T^*M)$ and $\psieps^{-2/(n-1)-1}\nabla_{e_j} \chi$ is bounded the integral involving  $\nabla_{e_j} \chi$ is well defined. For the integral involving $\nabla^M_{e_j}\psieps$ we observe 
 	\begin{align}
 		\begin{split}
 			&\sum_{\underset{j\neq k }{j,k=1}}^{n} \left(\nabla^M_{e_j}\psieps\right) \left\langle e_j \cl \psi, e_k \cl \nabla^{\Sigma}_{e_k} \psi \right\rangle\\
 			=& \sum_{j,k=1}^n \left(\nabla^M_{e_j}\psieps\right) \left\langle e_j \cl \psi, e_k \cl \nabla^{\Sigma}_{e_k} \psi \right\rangle - \sum_{j=1}^{n} \nabla^M_{e_j} \psieps \left\langle  \psi, \nabla^{\Sigma}_{e_j} \psi\right\rangle  \\
 			=& \left\langle \nabla^M \psieps \cl \psi, D\psi \right\rangle - g\left(\nabla^M \psieps ,\left\langle \psi, \nabla^{\Sigma} \psi  \right\rangle \right). \label{eq:step3help3local}
 		\end{split} 
 	\end{align}
 	Hence, the diamagnetic inequality $|\nabla^{\Sigma}\psi| \leq |\nabla^M |\psi||$, $|\psi|/\psieps\leq 1$ and $\psieps\nabla^M \psieps = |\psi| \nabla^M |\psi|$ imply
 	\begin{align*}
 		&\int_M\psieps^{-\frac{2}{n-1}-1} \left| \sum_{\underset{j\neq k }{j,k=1}}^{n} \left(\nabla^M_{e_j}\psieps\right) \left\langle e_j \cl \psi, e_k \cl \nabla^{\Sigma}_{e_k} \psi \right\rangle  \right| \chi dx\\
 		=&\int_M\psieps^{-\frac{2}{n-1}-1}\left[  \left\langle \nabla^M \psieps \cl \psi, D\psi \right\rangle - g\left(\nabla^M \psieps ,\left\langle \psi, \nabla^{\Sigma} \psi  \right\rangle \right) \right] \chi dx\\
 		\leq& C_{\chi,\varepsilon,n}\int_M \psieps^{-\frac{2}{n-1}} \left[ \left|\nabla^{\Sigma} \psi \right| \left| D\psi\right|+\left|\nabla^{\Sigma}\psi \right|^2\right] dx\\
 		\leq&\tilde{C}_{\chi,\varepsilon,n}\int_M \left|\nabla^{\Sigma} \psi\right|^2 dx <\infty.
 	\end{align*}
 	Thus, every integral in \eqref{eq:step3integralident} is well defined.\\ By density of $C^{\infty}_{c}(M,\Sigma)$ in $H^1(M,\Sigma)$, we prove \eqref{eq:step3integralident} for $\psi \in C^{\infty}_{c}(M,\Sigma)$. 
 	To see \eqref{eq:step3integralident}, we observe that $h_{k}:= f\langle \psi, e_k e_j \cl \nabla^{\Sigma}_{e_j}\psi\rangle $ is compactly supported on $U\subset M$. For the vector field $Y=h_k e_k$ on $M$ we have $\text{div}(Y)dx=\nabla^M_{e_k}h_kdx$. Thus, the divergence theorem implies for any $k$
 	\begin{align}
 		0=\int_M \nabla^M_{e_k}h_kdx. \label{eq:divergencethm}
 	\end{align}
 	Hence, \eqref{eq:divergencethm} and metric compatibility of $\langle \cdot, \cdot \rangle$ and $\nabla^{\Sigma}$ yield
 		\begin{align}
 		\begin{split}
 			&\int_{M} f \left\langle \nabla^\Sigma_{e_k} \psi, e_k e_j\cl \nabla^\Sigma_{e_j} \psi \right\rangle dx+\int_{M} f \left\langle \nabla^\Sigma_{e_j} \psi, e_j e_k\cl \nabla^\Sigma_{e_k} \psi \right\rangle dx\\=& -\int_M\left( \left(\nabla^M_{e_k}f\right) \left\langle \psi, e_k e_j \cl \nabla^\Sigma_{e_j} \psi \right\rangle + \left(\nabla^M_{e_j}f\right) \left\langle \psi, e_j e_k \cl \nabla^\Sigma_{e_k} \psi \right\rangle \right)dx\\
 			&- \int_{M} \left( f \left\langle \psi, \nabla^\Sigma_{e_k}\left(e_ke_j \cl \nabla^\Sigma_{e_j}\psi\right)\right\rangle dx + f \left\langle \psi, \nabla^\Sigma_{e_j}\left(e_je_k \cl \nabla^\Sigma_{e_k}\psi\right)\right\rangle\right) dx.
 		\end{split} \label{eq:step3help1}
 	\end{align}
 	For the second integral let $p \in M$ be the origin of the synchronous framing then
 	\begin{align*}
 		\nabla_{e_j}e_k \big|_{p}=0
 	\end{align*}
 	for all $j$ and $k$. Summing in \eqref{eq:step3help1} over $k<j$, the last integral yields the pointwise equation
 \begin{align}
 	\begin{split}
 		f\left\langle \psi, \sum^n_{k<j} \left(\nabla^\Sigma_{e_k} \left(e_k e_j \nabla^\Sigma_{e_j} \psi\right) - \nabla^\Sigma_{e_j}\left(e_k e_j \nabla^\Sigma_{e_k} \psi \right)\right)\right\rangle\big|_p
 		=f\left\langle\psi, \mathcal{K}\psi \right\rangle\big|_p,
 	\end{split} \label{eq:step3help2}
 \end{align}
 where the right hand side is independent on the choice of the synchronous frame, cf. \cite[Definition 3.7]{roe}. Summing over $k<j$ in \eqref{eq:step3help1} and using \eqref{eq:step3help2} we conclude \eqref{eq:step3integralident}. 
Moreover, as $\chi$ is supported in a nice coordinate patch we can integrate over the local expression
\begin{align}
	\begin{split}
		&\left|\nabla^\Sigma \psi\right|^2- \left|D\psi\right|^2 \\
		=&- \sum^n_{k<j}\left( \left\langle e_j \cl \nabla^\Sigma_{e_j} \psi, e_k \cl \nabla^\Sigma_{e_k} \psi \right\rangle+ \left\langle e_k \cl \nabla^\Sigma_{e_k} \psi, e_j \cl \nabla^\Sigma_{e_j} \psi \right\rangle \right),
		\end{split} \label{eq:step3help4}
	\end{align}
insert \eqref{eq:step3integralident} and obtain
\begin{align}
		\begin{split}
			&\int_M \left|\nabla^\Sigma \psi\right|^2f dx-\int_M \left|D\psi\right|^2 f dx \\
				=&\sum^n_{k<j} \int_{M} \left( f \left\langle \nabla^\Sigma_{e_j}\psi, e_je_k \cl \nabla^\Sigma_{e_k}\psi\right\rangle dx + f \left\langle \nabla^\Sigma_{e_k} \psi, e_ke_j \cl \nabla^\Sigma_{e_j}\psi\right\rangle\right) dx\\
			=&\sum^n_{\underset{j\neq k }{j,k=1}} \int_M \left(\nabla^M_{e_k}f\right)\left\langle e_k\cl \psi,e_j \cl \nabla^\Sigma_{e_j} \psi \right\rangle dx -\int_M f\left\langle \psi,\mathcal{K} \psi \right\rangle dx.
			\end{split} \label{eq:step3help5}
\end{align}
Inserting \eqref{eq:step3help3} into \eqref{eq:step3help5} yields
\begin{align}
	\begin{split}
		&\sum^n_{\underset{j\neq k }{j,k=1}} \int_M \left(\nabla^M_{e_k}f\right)\left\langle e_k\cl \psi,e_j \cl \nabla^\Sigma_{e_j} \psi \right\rangle dx -\int_M f\left\langle \psi,\mathcal{K} \psi \right\rangle dx  \\
		=&\sum^n_{\underset{k\neq j }{k,j=1}}\int_{M} \Big( \psieps^{-\frac{2}{n-1}} \left(\nabla^M_{e_j} \chi\right) \left\langle e_j \cl \psi, e_k \cl \nabla^\Sigma_{e_k} \psi \right\rangle\\ &- \frac{2}{n-1} \psieps^{-\frac{2}{n-1}-1} \left(\nabla^M_{e_j} \psieps\right) \left\langle e_j \cl \psi, e_k \cl \nabla^\Sigma_{e_k} \psi\right\rangle \Big) dx -\int_M f \left\langle \psi, \mathcal{K}\psi \right\rangle dx \label{eq:step3help6}
	\end{split}
\end{align}
We take real parts on both sides and observe that the first and last integral on the right hand side in \eqref{eq:step3help6} are exactly the two last integrals in \eqref{eq:step3}. For the integral involving $\nabla^M_{e_j} \psieps$, we insert \eqref{eq:step3help3local} into the integral and obtain 
\begin{align}
	\begin{split}
		-\frac{2}{n-1}&\int_M \sum_{\underset{j\neq k }{k,j=1}}^n \psieps^{-\frac{2}{n-1}-1}\left(\nabla \psieps\right) \left\langle e_j \cl \psi, e_k \cl \nabla^{\Sigma}_{e_k} \psi \right\rangle  \chi dx\\
		=-\frac{2}{n-1}&\int_M \psieps^{-\frac{2}{n-1}-1} \text{Re} \left(\left\langle (\nabla^M \psieps),D\psi \right\rangle \right) \chi dx \\
		+\frac{2}{n-1}&\int_M \psieps^{-\frac{2}{n-1}-1} g\left((\nabla^M\psieps), \text{Re}\left(\left\langle \psi, \nabla^{\Sigma}\psi \right\rangle\right)\right)\chi dx.
	\end{split} \label{eq:step3helpnew}
\end{align}
We observe that the first integral in \eqref{eq:step3helpnew} is exactly the third term on the right hand side of \eqref{eq:step3}.
Imposing \eqref{eq:step0} in the last integral in \eqref{eq:step3helpnew} yields
\begin{align*}
	\psieps^{-\frac{2}{n-1}-1}g\left(\nabla^M\psieps,\text{Re}\left(\left\langle \psi, \nabla^{\Sigma} \psi  \right\rangle\right) \right)=\psieps^{-\frac{2}{n-1}} \left|\nabla^M \psieps\right|^2.
\end{align*}
To finish the proof of step 3 we note that rewriting the right hand side yields
\begin{align*}
	\frac{2}{n-1}\psieps^{-\frac{2}{n-1}} \left| \nabla^M\psieps\right|^2 =\frac{2(n-1)}{(n-2)^2}\left| \nabla^{M}\left( \psieps^{\frac{n-2}{n-1}}\right)\right|^2
\end{align*}
and \eqref{eq:step3} follows. \\
For step (4) we claim
\begin{align}
	\begin{split}
		\int_{M} \psieps^{-\frac{2}{n-1}}\left|\nabla^\Sigma\psi\right|^2dx =& \frac{2(n-1)}{(n-2)^2} \int_{M} \left|\nabla^M\left( \psieps^{\frac{n-2}{n-1}}\right) \right|^2dx \\
		+&\int_{M}\psieps^{-\frac{2}{n-1}} \left|D\psi\right|^2 dx \\
		-& \frac{2}{n-1}\int_{M} \psieps^{-\frac{2}{n-1}-1} \text{Re}\left( \left\langle \left(\nabla^M \psieps \right)\cl \psi, D\psi \right\rangle \right)dx \\
		-& \int_{M}\psieps^{-\frac{2}{n-1}} \text{Re}\left( \left\langle \psi, \mathcal{K}\psi \right\rangle \right).
	\end{split} \label{eq:step4}
\end{align}
To see \eqref{eq:step4} we cover $M$ by sets $U_l$ with the properties in step (3) and choose a partition of unity $\chi_l$ subordinate to $U_l$. We apply \eqref{eq:step3} for $\chi_l$ and sum over $l$. We note that
\begin{align*}
	\sum_l \sum_{\underset{j\neq k}{j,k=1}}^{n}&\int_{M} \psieps^{-\frac{2}{n-1}}\text{Re}\left(\left\langle e_j \cl \psi, e_k \cl \nabla^\Sigma_{e_k} \psi\right\rangle\right) \left(\nabla^M_{e_j} \chi_l\right) dx\\
	= \sum_{\underset{j\neq k}{j,k=1}}^{n}&\int_{M} \psieps^{-\frac{2}{n-1}}\text{Re}\left(\left\langle e_j \cl \psi, e_k \cl \nabla^\Sigma_{e_k} \psi\right\rangle\right) \left(\nabla^M_{e_j}\sum_{l} \chi_l\right) dx=0,
\end{align*}
as $\sum_l \chi_l \equiv 1$.
Collecting terms we arrive at \eqref{eq:step4}.\\ To finish we use \eqref{eq:step1} and \eqref{eq:step2} and observe
\begin{align*}
	&\left|\nabla^{\Sigma} \left( \frac{\psi}{\psieps^{\frac{n}{n-1}}}\right)\right| \psieps^2 - \frac{1}{n} \left|-iD \left(\frac{\psi}{\psieps^{n/n-1}}\right)\right|^2\psieps^2 \\
	=&\frac{\left|\nabla^{\Sigma} \psi\right|^2}{\psieps^{\frac{2}{n-1}}} + \left| \nabla^M \left( \psieps^{\frac{n-2}{n-1}} \right) \right|^2\left[ \left(\frac{n}{n-2}\right)^2\frac{\left|\psi\right|^2}{\psieps^2} - \frac{2n(n-1)}{(n-2)^2} \right] \\
	+&\frac{\left|D\psi\right|^2}{\psieps^{\frac{2}{n-1}}} -\frac{1}{n} \left(\frac{n}{n-2}\right)^2\left|\nabla^M \left(\psieps^{\frac{n-2}{n-1}}\right)\right|^2\frac{\left|\psi\right|^2}{\psieps^2} \\ +&\frac{2}{n-1}\frac{1}{\psieps^{\frac{2}{n-1}+1}}\text{Re} \left( \left\langle D\psi, \nabla^M \psieps \cl \psi \right\rangle \right).
\end{align*}
Integrating over $M$, using \eqref{eq:step4} to express the first and last term above in terms of $|D\psi|^2$ and $|\nabla^M\big(\psieps^{\frac{n-2}{n-1}}\big)|^2$, here also the term involving $\mathcal{K}$ appears, and rewriting the identity yields \eqref{eq:integraliden}.
\end{proof}
\rem[]{Proposition \ref{prop:InteIdent} is still true if we exchange the assertion of $M$ being closed by $M$ being complete.}
\rem[]{
 	We emphasize that the conditions in \cite[Proposition 7]{FrankLoss2} avoid the assertion $\psi \in L^{2}(M,\Sigma)$ and it is already sufficient to suppose $\nabla \psi \in D^{1}(M,\Sigma)$, the space of all functions where $\nabla \psi \in L^2(M,T^*M\otimes \Sigma)$, $\psi \in L^1_{\text{loc}}(M,\Sigma)$ and $|\{|\psi|>\tau\}|<\infty$ for all $\tau>0$. This is based on the fact that $(\mathbb{R}^n,g_{\mathbb{R}^n})$ satisfies a Hardy inequality and one can cover $(\mathbb{R}^n,g_{\mathbb{R}^n})$ by one global chart. 
 	Thus, in step 4 in the proof on $\mathbb{R}^n$, it is possible scale the cut-off function $\chi$ by some factor $R$ such that it is exhausting. 
 	Then, after imposing certain bounds on the gradient of the cut-off function, one applies the Hardy inequality and dominated convergence yields the same result. However,
 	as discussed in Remark \ref{rem:l2assum}, the assumption $\psi \in H^{1}(M,\Sigma)$ does not affect the generality of the assumptions in our main theorem.}
	\section{\textsc{Main theorem}\label{sec:mainthm}}
	\subsection{The compact case}
	In this section, we prove Theorem \ref{thm:spin}. In Theorem \ref{thm:mainyam} we first prove a slightly more general version of Theorem \ref{thm:spin}. By the means of Remark \ref{rem:spin}, Theorem \ref{thm:spin} follows as a corollary. 
	\thm[]{Let $(M,g)$ be a closed Riemannian manifold of dimension $n\geq3$. Suppose $\psi \in L^{p}(M,\Sigma)$ for some $p$ with $\frac{2n}{n-2}\leq p <\infty$ is a nontrivial solution to $D\psi =i A\cl \psi$,
		then 
		\begin{align}
			\frac{n}{4(n-1)}Y(M,[g]) \leq \big\|  |A|^2 \text{id} - \frac{n}{n-1}F^{\Sigma}\big\|_{L^{\frac{n}{2}}}. \label{eq:inequYam}
		\end{align} \label{thm:mainyam}
	}\\ \upshape
	\\ The proof follows the same strategy as the proof of \cite[Theorem 1]{FrankLoss2}. Since employing the Yamabe-invariant yields some technical issues, we will include the full proof here for the sake of better readability. 
		\begin{proof}[Proof of Theorem \ref{thm:mainyam}]
			By compactness of $M$ and $\psi \in L^{p}(M,\Sigma)$ for some $p \in [\frac{2n}{n-2},\infty)$ we obtain due to Hölder's inequality
			\begin{align*}
				\psi \in L^{\tilde{p}}(M,\Sigma) \quad \text{for all } \tilde{p}  \in [2,p].
			\end{align*}
			Hence, \eqref{eq:maineq}, $A\in L^n(M,TM)$ and Hölder's inequality imply $D\psi \in L^2(M,\Sigma)$ and Garding's inequality, cf. Proposition \ref{rem:gard}, yields $\psi \in H^{1}(M,\Sigma)$.
			By \eqref{eq:twistor_abs}, we rewrite the left hand side in our integral identity and observe
			\begin{align}
				|T_g \varphi|^2=  \sum_{j=1}^{n} \Big|\big[\nabla_{e_j} + \frac{1}{n} e_j \cl D  \big]\varphi\Big|^2\geq 0  \label{eq:pmainthm1}
			\end{align}
			for $\varphi \in \Gamma(M,\Sigma)$. 
			Making use of \eqref{eq:pmainthm1} and $|\psi|\leq \psieps$, we rearrange \eqref{eq:integraliden} and obtain
			\begin{align*}
				\frac{n-1}{n-2} \int_M \left| \nabla^{M} \left(\psieps^{\frac{n-2}{n-1}}\right) \right|^2 dx\leq \frac{n-1}{n} \int_M \frac{|D\psi|^2}{\psieps^{\frac{2}{n-1}}}dx - \int_M \frac{\left\langle\mathcal{K} \psi,\psi \right\rangle}{\psieps^{\frac{2}{n-1}}}dx 
			\end{align*}
			After inserting \eqref{eq:maineq} and $\mathcal{K}=F^{\Sigma} +\frac{\text{scal}_g}{4}$ the inequality above can be rewritten as
			\begin{align}
				\begin{split}
					\frac{n-1}{n-2} \int_M \left| \nabla^{M} \left(\psieps^{\frac{n-2}{n-1}}\right) \right|^2 dx\leq& \frac{n-1}{n} \int_M \frac{\left\langle \left(|A|^2\text{id} -\frac{n}{n-1} F^{\Sigma}\right)\psi,\psi\right\rangle}{\psieps^{\frac{2}{n-1}}}dx \\
					&-\int_M \frac{\text{scal}_g}{4} \frac{|\psi|^2}{\psieps^{\frac{2}{n-1}}}dx
				\end{split} \label{eq:pmainthm2}
			\end{align}
			We define the function 	$\widetilde{|\psi|}_{\varepsilon} : = \psieps^{\frac{n-2}{n-1}}- \varepsilon^{\frac{n-2}{n-1}}$ and observe $\nabla^M \widetilde{|\psi|}_{\varepsilon} = \nabla^M \psieps^{\frac{n-2}{n-1}}$. 
			After adding $\int_M \frac{\text{scal}_g}{4}\widetilde{|\psi|}^2_{\varepsilon}dx$ on both sides in \eqref{eq:pmainthm2}, we estimate the left hand side in \eqref{eq:pmainthm2} by inserting the definition of the Yamabe invariant. This yields
			\begin{align*}
				\frac{Y(M,[g])}{4}\big(\int_{M}\widetilde{|\psi|}^{\frac{2n}{n-2}}_{\varepsilon}\big)^{\frac{n-2}{n}} 
				\leq&
				\frac{n-1}{n}\int_M \frac{ \langle\left(|A|^2\text{id}- \frac{n}{n-1}F^\Sigma\right)\psi,\psi\rangle}{\psieps^{\frac{2}{n-1}}} \\
				&+\int_M \frac{\text{scal}_g}{4} \bigg[\widetilde{|\psi|}^2_{\varepsilon}-\frac{|\psi|^2}{\psieps^{\frac{2}{n-1}}} \bigg]
			\end{align*}
		It remains to understand the behavior of $\varepsilon\rightarrow 0$ on both sides. We have three estimates
		\begin{enumerate}
			\item[(i)] as $\psieps$ is non-increasing in $\varepsilon$, by monotone convergence we have 
			\begin{align*}
				\lim_{\varepsilon\rightarrow 0 } \int_M \widetilde{|\psi|}_{\varepsilon}^{\frac{2n}{n-2}} dx= \int_M |\psi|^{\frac{2n}{n-1}} dx,
			\end{align*}
			\item[(ii)] Hölder's inequality and $\psieps^{-\frac{2}{n-1}}|\psi|^2\leq |\psi|^{\frac{2(n-2)}{n-1}}$ imply
			\begin{align*}
				\int_M \frac{\langle\left(|A|^2\text{id} - \frac{n}{n-1}F^\Sigma\right)\psi,\psi\rangle}{\psieps^{\frac{2}{n-1}}} \leq \big\|  |A|^2\text{id} - \frac{n}{n-1}F^\Sigma \big\|_{L^{\frac{n}{2}}} \big(\int_M |\psi|^{\frac{2n}{n-1}}dx\big)^{\frac{n-2}{n}},
			\end{align*}
			\item[(iii)] $\psi \in L^{\frac{2(n-2)}{n-1}}(M,\Sigma)$, dominated convergence and $n\geq 3$ yield
			\begin{align*}
				\lim_{\varepsilon\rightarrow 0} 	\int_M \frac{\psieps^2-|\psi|^2}{\psieps^{\frac{2}{n-1}}} dx=  \int_M \lim_{\varepsilon\rightarrow 0}\frac{\varepsilon^2}{(|\psi|^2+ \varepsilon^2)^{\frac{1}{n-1}}}dx=0. 
			\end{align*}
		\end{enumerate}
		Thus, the assertion that $\psi$ is nontrivial implies \eqref{eq:inequYam}.
		\end{proof}
	\rem[]{In the proof we directly see that imposing $\psi \in L^{p}(M,\Sigma)$ for $p\in [2n/(n-2),\infty)$ directly implies $\psi \in L^2(M,\Sigma)$. Thus, for our purpose the assumption $\psi \in L^2(M,\Sigma)$ in Proposition \ref{prop:InteIdent} does not affect the generality of our result.\label{rem:l2assum}}
	\rem[]{Compared to the assumption in \cite[Theorem 1]{FrankLoss2}, we have to replace the condition $p\in(n/(n-1),\infty)$ with $p\in [2n/(n-2),\infty)$. On $\mathbb{R}^n$ in \cite[Theorem 2.1]{FrankLoss} the authors proved that $\psi\in L^{p}(\mathbb{R}^n,\mathbb{C}^N)$ for some $p \in (n/(n-1),\infty)$, being a nontrivial solution to \eqref{eq:maineq}, already implies $\psi\in L^{p}(\mathbb{R}^n,\mathbb{C}^N)$ for all $p \in (n/(n-1),\infty)$. The proof of this theorem is based on certain techniques which are valid on $\mathbb{R}^n$ but do not provide very general analogues on manifolds, such as the Hardy-Littlewood-Sobolev inequality and explicit formulas for the Green function of the Dirac operator. 
	The question whether a statement similar to \cite[Theorem 2.1]{FrankLoss} holds on (non-)compact Riemannian manifolds has yet not been answered.} 
	\subsection{The non-compact case}
	Imposing slightly stronger assumptions as in Theorem \ref{thm:mainyam}, we can prove a similar theorem for manifolds of bounded geometry. As before Theorem \ref{thm:mainboundedgeomspin} follows as a corollary by the means of Remark \ref{rem:spin}.
	\thm[]{Let $(M,g)$ be a Riemannian manifold of bounded geometry, $\dim M=n\geq3$, and $(\Sigma,\langle \cdot,\cdot\rangle,\nabla^{\Sigma})$ a Clifford bundle of bounded geometry over $(M,g)$ such that $F^{\Sigma}\in L^{\frac{n}{2}}(M,End(S))$. Assume $A\in L^{n}(M,TM)$ and $\|A\|_{L^{\infty}}<\infty$. If $\psi \in L^{2}(M,\Sigma)$ is a nontrivial solution of $D\psi = i A\cl \psi$, then
	\begin{align*}
		\frac{n}{4(n-1)}Y(M,[g]) \leq  \big\|  |A|^2\text{id}- \frac{n}{n-1}F^{\Sigma}  \big\|_{L^{\frac{n}{2}}} .
		\end{align*} \label{thm:mainboundedgeom}}
	\begin{proof}
		The assertions $\|A\|_{L^{\infty}}<\infty$ and $\psi \in L^2(M,\Sigma)$ imply $\|D\psi\|_{L^2}<\infty$ and hence, by \cite[Lemma A.1.4]{shu} and Proposition \ref{rem:gard}, $\psi \in H^{1,2}(M,\Sigma)$. By \cite[Theorem 21]{gro2} $H^{1,2}(M,\Sigma)\subset L^{2n/(n-2)}(M,\Sigma)$ holds for manifolds and vector bundles of bounded geometry. Thus, by interpolation of Sobolev spaces, we have $\psi \in L^{p}(M,\Sigma)$ for $p \in \{2n/(n-2),2n/(n-1),2(n-2)/(n-1),2\}$. The rest of the proof now follows exactly the proof of Theorem \ref{thm:mainyam}.
	\end{proof}
		\ex[]{One important class of manifolds which is covered by Theorem \ref{thm:mainboundedgeom} are manifolds conformally equivalent to the hyperbolic space $\mathbb{H}^n$ equipped with the standard metric $g_{\mathbb{H}^n}$. By \cite[Lemma 7.1]{amm1} the Yamabe invariant of $(\mathbb{H}^n,g_{\mathbb{H}^n})$ coincides with the Yamabe invariant $(\mathbb{S}^n,g_{\mathbb{S}^n})$ and we observe that the constants appearing in the inequality are the same for the classical model spaces. \\
		A second class of manifolds for which Theorem \ref{thm:mainboundedgeom} is applicable are asymptotically flat (ALE) manifolds. ALE manifolds almost look like the Euclidean space outside of a compact set and are a central subject studied in many current problems in geometry and analysis.}
	\section{\textsc{The  case of equality}\label{sec:recovery}}
	In this section, we will first recover the case of equality in \eqref{eq:inequYamspin} on $(\mathbb{S}^n,g_{\mathbb{S}^n})$ for odd $n$. Afterwards, we give a  characterization of closed manifolds on which equality in \eqref{eq:inequYamspin} is possibly attained. Throughout we assume that $(M,g,\sigma)$ is a closed Riemannian spin manifold with fixed spin structure $\sigma$.\\
	We first recall some facts about the case of equality on $\mathbb{R}^n$ from \cite[Section 1]{FrankLoss2}.
	\rem[]{\label{rem:equality_on_rn}Let $n$ be odd.
	According to \cite[Theorem 5]{FrankLoss2} the pair $(A,\psi)$ solves \eqref{eq:maineq} on $\mathbb{R}^n$ and attains equality in \eqref{eq:inequYamspin} for 
	\begin{align*}
		A(x)=n\left(\frac{1}{1+|x|^2}\right)^2 \left(\left(1-|x|^2\right)e_1+ 2g_{\mathbb{R}^n}( x,e_1)x + 2 L_S x \right) 
	\end{align*} 
	and
	\begin{align*}
		\psi(x)= \left(\frac{1}{1+|x|^2}\right)^{\frac{n}{2}} \left(1+ i s x \cl \right)\psi_0,
	\end{align*}
	where, $s \in \left\{1,-1
	\right\}$,  $e_1$ denotes the standard canonical first basis vector on $\mathbb{R}^n$ and $L_S$ is the skew-symmetric matrix
	\begin{align*}
	L_S=\begin{bmatrix}
		 0 &  &  &  &  &  &  \\
		     & 0 & -1 &  &  &  &  \\
		     & 1 & 0 &  &  &  &  \\
		     &  &  & \ddots &  & &  \\
		     &  &  &  &  &0 & -1\\
		     &  &  & &  &  1 & 0 \\ 
	\end{bmatrix}.
	\end{align*}
	Here $\psi_0$ a constant spinor such that $|\psi_0|=1$.}
	\rem[The case $(\mathbb{S}^n,g_{\mathbb{S}^n})$]{Let $n$ be odd. 
	We first
	note that $(\mathbb{S}^n\setminus \{P\},g_{\mathbb{S}^n})$ is conformally equivalent to $(\mathbb{R}^n,g_{\mathbb{R}^n})$ for some $P\in \mathbb{S}^n$ via stereographic projection. Without loss of generality, we choose $P=N$ as the northpole of $\mathbb{S}^n$. \\
	There is an isomorphism of spinor  bundles  $F: \Sigma(\mathbb{R}^n,g_{\mathbb{R}^n},\sigma)\rightarrow \Sigma(\mathbb{S}^n\setminus\{N\},g_{\mathbb{S}^n},\sigma)$, see e.g. \cite[Proposition 2.2.1]{amm2}. Then, the conformal transformation rule, see e.g. \cite[Proposition 2.2.1]{amm2}, for the Dirac operator yields
	\begin{align}
		D_{g_{\mathbb{S}^n}}F\left(\psi\right) = F\left(\left(\frac{2}{1+|x|^2}\right)^{-\frac{n+1}{2}} D_{g_{\mathbb{R}^d}}\left(\left( \frac{2}{1+|x|^2}\right)^{\frac{n-1}{2}} \psi \right)\right).\label{eq:confDir}
	\end{align}
	We define $\varphi:= F\left(\left(\frac{2}{1+|x|^2}\right)^{-\frac{n-1}{2}}\psi\right)$.  Then, \eqref{eq:maineq} and \eqref{eq:confDir} imply
	\begin{align*}
		D_{g_{\mathbb{S}^n}}\varphi= \frac{1+|x|^2}{2} A \cl \varphi. 
	\end{align*}
	Thus, the pair $\left(\bar{A}:=\frac{1+|x|^2}{2} A ,\varphi\right)$ solves \eqref{eq:maineq} on $\mathbb{S}^n\setminus\{N\}$. We emphasize that the only condition on $\bar{A}$ and $\varphi$ in Theorem \ref{thm:spin} is $\bar{A}\in L^n(M,TM)$ and $\varphi \in L^p(M,\Sigma)$ for some $p$ large enough. Hence, as the $L^n$ norm of $A$ is conformally invariant, setting $\bar{A}$ to zero at $N$ yields 
	\begin{align*}
		\|\bar{A}\|^2_{L^n}= \|A\|^2_{L^n}= \frac{n}{4(n-1)}Y(\mathbb{S}^n,[g_{\mathbb{S}^n}]).
	\end{align*}
Hence, we have recovered equality on odd dimensional spheres. \label{rem:recovery_on_sphere}}
\rem[]{By a straightforward computation based on the conformal transformation rule we can recover the case of equality on $\mathbb{H}^n$ from $\mathbb{R}^n$, as by \cite[Lemma 7.1]{amm1} $Y(\mathbb{H}^n,[g_{\mathbb{H}^n}])=Y(\mathbb{S}^n,[g_{\mathbb{S}^n}])=\frac{4(n-1)}{n-2}S_n$. }
\\
We will now give a characterization of closed manifolds where equality is possibly attained in \eqref{eq:inequYamspin}. The proof follows the same strategy as the proof of \cite[Proposition 8]{FrankLoss2}, we will therefore only present a sketch of the proof. 
\thm[]{\label{thm:equality}Let $(M,g,\sigma)$ be a closed  Riemannian spin manifold of dimension $n\geq3$. Assume $A\in L^{n}(M,TM)$ and $0 \neq \psi \in L^{p}(M,\Sigma)$ is a solution to \eqref{eq:maineq}
	for some $p\in [2n/(n-2),\infty)$.
	If
	\begin{align}
		\|A\|^2_{L^n}=\frac{n}{4(n-1)}Y(M,[g]), \label{eq:equality}
	\end{align}
	then $|\psi|^{\frac{n-2}{n-1}}$ is an optimizer for the Yamabe functional, i.e., $|\psi|$ is smooth, $|\psi|>0$  and $\tilde{g}=|\psi|^{\frac{4}{n-1}}g$ is a complete metric with constant scalar curvature $Y(M,[g])$. Moreover $\psi/|\psi|^{n/(n-1)}$ belongs to the kernel of $T_g$.
}
\begin{proof}
		We first assume that $|\psi|> 0$ and that $|\psi|$ is smooth. At the end we will spend some words about the possible zeros of $|\psi|$. We define
		\begin{align*}
			P:= \int_M \left|T_g \left(\frac{\psi}{|\psi|^{\frac{n}{n-1}}}\right)\right|^2 \left|\psi\right|^2 dx.
		\end{align*}
		As $|\psi|>0$ we can set $\varepsilon=0$ in the integral identity \eqref{eq:integraliden} and obtain 
		\begin{align}
			P=\frac{n-1}{n}\int_M \left| A \right|^2 \left|\psi \right|^{\frac{2(n-2)}{n-1}}dx&-\frac{n-1}{n-2}\int_M\left| \nabla^{M}\left(|\psi|^{\frac{n-2}{n-1}}\right) \right|^2 dx\\
			&-\int_M \frac{\text{scal}_g}{4} |\psi|^{\frac{2(n-2)}{n-1}}dx. \label{eq:thm_equality_1}
		\end{align}
		Rearranging \eqref{eq:thm_equality_1} in terms of Hölder and Yamabe inequality yields
		\begin{align*}
			P+R^{(1)}+R^{(2)}=S,
		\end{align*}
		where
		\begin{align*}
			R^{(1)}&=\frac{n-1}{n} \left[\|A\|^2_{L^n} \left(\int_M |\psi|^{\frac{2n}{n-1}}dx \right)^{\frac{n-2}{n}}- \int_M |A|^2 |\psi|^{\frac{2(n-2)}{n-1}} dx\right], \\
			R^{(2)}&=\frac{n-1}{n-2}\bigg[\int_M \left|\nabla^{M} \left(|\psi|^{\frac{n-2}{n-1}}dx\right) \right|^2 +\frac{1}{c_n}\int_M \text{scal}_g |\psi|^{\frac{2(n-2)}{n-1}}dx\\
			&- Y(M,[g]) \left(\int_M |\psi|^{\frac{2n}{n-1}}dx\right)^{\frac{n-2}{n}}\bigg]
		\end{align*}
		and
		\begin{align*}
			S=\left( \frac{n-1}{n} \|A\|^2_{L^n}-\frac{Y(M,[g])}{4}\right)\left(\int_M |\psi|^{\frac{2n}{n-1}}dx\right)^{\frac{n-2}{n}}
		\end{align*}
		for $c_n=\frac{4(n-1)}{n-2}$. According to the definition of the Yamabe invariant, Hölder's inequality and Remark \ref{rem:twistor}, we have $R^{(1)},R^{(2)},P\geq0$ and hence $S\geq0$. Thus,
			\begin{align*}
			S= \left(\frac{n-1}{n} \|A\|^2_{L^n}-\frac{Y(M,[g])}{4}\right)\left(\int_M |\psi|^{\frac{2n}{n-1}}dx\right)^{\frac{n-2}{n}}\geq0
		\end{align*}
		and as $\psi \neq0 $ we conclude \eqref{eq:inequYamspin} again. As we assume $\|A\|^2_{L^n}=\frac{n}{4(n-1)}Y(M,[g])$, we obtain $S=0$ and hence $P=R^{(1)}=R^{(2)}=0$. For $R^{(1)}=0$, we need equality in Hölder's inequality and hence,
		\begin{align*}
			|A|=\lambda |\psi|^{\frac{2}{n-1}}
		\end{align*}
		for some constant $\lambda$. If $R^{(2)}=0$, then due to the definition of the Yamabe invariant $|\psi|^{\frac{n-2}{n-1}}$ is an optimizer for the Yamabe problem  and hence $(M,\tilde{g})$, where $\tilde{g}= |\psi|^{\frac{4}{n-1}}g$, has constant scalar curvature $\text{scal}_{\tilde{g}}=Y(M,[g])$. Finally, as $P=0$, we have $T_g\left(\frac{\psi}{|\psi|^{\frac{n}{n-1}}}\right)=0$. \\
		Let us show that $|\psi|$ indeed does not have any zeros. 
		We again use the regularization $\psieps$ introduced in Section \ref{sec:intident} and rewrite the integral identity in a similar manner as above. Now $P=P_{\varepsilon},R^{(1)}=\repsone,R^{(2)}=\repstwo$ and $S=S_{\varepsilon}$ are all dependent on $\varepsilon$ and an additional term $R_{\varepsilon}$ appears, where
			\begin{align*}
			R_{\varepsilon}= \frac{n(n-1)}{(n-2)^2} \int_M \left|\nabla^M\left( \psieps^{\frac{n-2}{n-1}}\right)\right|^2 \frac{\varepsilon^2}{\psieps^2}dx.
		\end{align*}
		We directly observe that $R_{\varepsilon}\geq0$. Moreover, we write $S_{\varepsilon}=S^{(1)}_{\varepsilon}+S^{(2)}_{\varepsilon}$, given by
			\begin{align*}
			S^{(1)}_{\varepsilon}=&\frac{n-1}{n} \|A\|^2_{L^n}\left(\int_M |\psi|^{\frac{2n}{n-2}}\psieps^{-\frac{2n}{(n-1)(n-2)}}dx
			\right)^{\frac{n-2}{n}}\\
			&-\frac{n-1}{n-2}\frac{Y(M,[g])}{c_n} \left(\int_M \widetilde{|\psi|}^{\frac{2n}{n-2}}_{\varepsilon}dx
			\right)^{\frac{n-2}{n}},\\
			S^{(2)}_{\varepsilon}=& \int_M \frac{\text{scal}_g}{4} \left(
			\frac{|\psi|^2}{\psieps^{\frac{2}{n-1}}} -\widetilde{|\psi|}^2_{\varepsilon}dx
			\right).
		\end{align*}
		Hence, as $R_\varepsilon,\peps,\repsone,\repstwo\geq0$, we have
		\begin{align*}
			S^{(1)}_{\varepsilon}\geq -S^{(2)}_{\varepsilon}
		\end{align*}
		and by monotone convergence and $\psi \in L^{\frac{2n}{n-1}}(M,\Sigma)$ we conclude
			\begin{align*}
			\left(
			\frac{n-1}{n-2}\left(\|A\|^2_{L^n}- \frac{Y(M,[g])}{c_n}\right) \left(\int_M|\psi|^{\frac{2n}{n-1}}dx \right)^{\frac{n-2}{n}}
			\right) &=\lim_{\varepsilon\rightarrow 0}S^{(1)}_{\varepsilon}  \\
			&\geq \lim_{\varepsilon\rightarrow 0}-S^{(2)}_{\varepsilon}=0.
		\end{align*}
		As we assume that $\psi$ is nontrivial, we obtain
		\begin{align*}
			\|A\|^2_{L^n}\geq\frac{n}{4(n-1)}Y(M,[g]).
		\end{align*}
		Imposing \eqref{eq:equality} yields
		\begin{align*}
			\lim_{\varepsilon\rightarrow 0}R_{\varepsilon}=\lim_{\varepsilon\rightarrow0}P_{\varepsilon}=\lim_{\varepsilon\rightarrow0}R^{(1)}_{\varepsilon}
			=\lim_{\varepsilon\rightarrow0}R^{(2)}_{\varepsilon}=0.
		\end{align*}
		For our purpose particularly interesting is the term $\repstwo$ given by
		\begin{align*}
			R^{(2)}_{\varepsilon}&= \frac{n-1}{n-2} \Big[
			\int_M\left( \left|\nabla^{M}\psieps^{\frac{n-2}{n-1}}\right|^2 +\frac{1}{c_n} \text{scal}_g\widetilde{|\psi|}_\varepsilon \right)dx\\
			&-\frac{Y(M,[g])}{c_n}\left(\int_M \widetilde{|\psi|}_{\varepsilon}^{\frac{2n}{n-2}}\right)^{\frac{n-2}{n}}\Big],
		\end{align*}
		where $\widetilde{|\psi|}_{\varepsilon}=|\psi|^{\frac{n-2}{n-1}}-\varepsilon^{\frac{n-2}{n-1}}$.
		By monotone convergence
		\begin{align*}
			\left(\int_{M}\widetilde{|\psi|}^{\frac{2n}{n-2}}_{\varepsilon} \right)^{\frac{n-2}{n}}\underset{\varepsilon\rightarrow0}{\rightarrow}
			\left(\int_M |\psi|^{\frac{2n}{n-1}}\right)^{\frac{n-2}{n}}
		\end{align*}
		and 
		\begin{align*}
			\frac{1}{c_n}\int_M\text{scal}_g \widetilde{|\psi|}^2_{\varepsilon} \underset{\varepsilon\rightarrow0}{\rightarrow} \frac{1}{c_n}\int_M \text{scal}_g |\psi|^{\frac{2(n-2)}{n-1}}.
		\end{align*}
			The fact that $\repstwo\rightarrow 0$ implies in particular that $\repstwo$ remains bounded. Thus,
		$\int_M \left|\nabla^{M} \psieps^{\frac{n-2}{n-1}}\right|^2dx$ remains bounded. Thus, as for $\varepsilon\rightarrow0$,
		we have $\psieps^{\frac{n-2}{n-1}} \rightarrow |\psi|^{\frac{n-2}{n-1}}$ in $L^{1}_{\text{loc}}$,
		we observe that $|\psi|^{\frac{n-2}{n-1}}$ is weakly differentiable. By \cite[Lemma 9]{FrankLoss2},
		we obtain
		\begin{align*}
			\int_M \left| \nabla^{M} |\psi|^{\frac{n-2}{n-1}}\right|^2dx 
			\leq
			\liminf_{\varepsilon\rightarrow 0}
			\int_M |\nabla^M \psieps^{\frac{n-2}{n-1}}|^2dx.
		\end{align*}
		Thus,
		\begin{align*}
			\underset{\varepsilon\rightarrow 0}{\liminf} \repstwo \geq& 
			\int_M |\nabla^M |\psi|^{\frac{n-2}{n-1}}|^2dx 
			+ \frac{1}{c_n} \int_M \text{scal}_g |\psi|^{\frac{2(n-2)}{n-1}} dx\\
			&- \frac{Y(M,[g])}{c_n} \left(\int_M |\psi|^{\frac{2n}{n-1}} dx\right)^{\frac{n-2}{n}}.
		\end{align*}
		As $\lim_{\varepsilon\rightarrow 0}\repstwo=0$, we obtain
		\begin{align*}
			\frac{Y(M,[g])}{c_n} \left(
			\int_M |\psi|^{\frac{2n}{n-1}}dx
			\right)^{\frac{n-2}{n}}=\int_M |\nabla^M |\psi|^{\frac{n-2}{n-1}}|^2dx +\frac{1}{c_n}\int_M \text{scal}_g|\psi|^{\frac{2(n-2)}{n-1}}dx.
		\end{align*}
		Hence, $|\psi|^{\frac{n-2}{n-1}}$ is an optimizer for the Yamabe problem. Thus, $|\psi|^\frac{n-2}{n-1}>0$. 
\end{proof}
\rem[]{Comparing the results in Theorem \ref{thm:equality} with the pair $(\bar{A},\varphi)$ that solves \eqref{eq:maineq} on $\mathbb{S}^n$, we note that $\varphi$ may not be extended continuously to all of $\mathbb{S}^n$ from $\mathbb{R}^n$, but after setting $\varphi(N)=F(2^{-\frac{n-1}{2}}\psi(0))$, the function $|\varphi|$ is defined on all of $\mathbb{S}^n$ and we have
\begin{align*}
	|\varphi|(x)= 2^{-\frac{n-1}{2}}
	\end{align*}
and indeed inserting $|\varphi|$ into the definition of the Yamabe invariant we obtain 
\begin{align*}
	\frac{Y(\mathbb{S}^n,[g_{\mathbb{S}^n}])}{c_n}\left(\int_{\mathbb{S}^n} |\varphi|^{\frac{2n}{n-1}}\right)^{\frac{n-2}{n}}= \int_{\mathbb{S}^n} |\nabla^{\mathbb{S}^n} |\varphi|^{\frac{n-2}{n-1}}|^2dx +\frac{1}{c_n}\int_{\mathbb{S}^n} \text{scal}_{g_{\mathbb{S}^n}}|\varphi|^{\frac{2(n-2)}{n-1}}dx.
\end{align*}
}Before we characterize the manifolds on which equality in \eqref{eq:inequYamspin} is possibly attained, we introduce Killing spinors.
\defi[{\cite[Section 5.2]{friedrich}}]{A spinor $\psi$ is called (real) Killing spinor if there exists a number $\mu \in \mathbb{C}$ ($\mathbb{R}$ for real Killing spinors), such that for all $X\in \mathfrak{X}(M)$
	\begin{align*}
		\nabla^{\Sigma}_X \psi =\mu X\cl \psi
	\end{align*}
	holds.}	
\rem[]{\label{rem:characterization}In this remark, we enlist classes of manifolds on which equality in \eqref{eq:inequYamspin} is possibly attained and observe that in $3$ dimensions $(M,g)$ has to be conformally equivalent to $(\mathbb{S}^3,g_{{\mathbb{S}^3}})$.\\ 
Recall from from the proof of Theorem \ref{thm:equality} that $\psi/|\psi|^{\frac{n}{n-1}} \in \text{ker}(T_g)$ and thus $M$ carries a nontrivial twistor spinor.
By \cite[Section 5.3]{friedrich}, after a possible conformal transformation, $\text{ker}(T_g)$ can be identified with the space of real Killing spinors on $(M,\bar{g})$. In \cite[Theorem 2', Theorem 3', Theorem 4', Theorem 5']{baer}
the author classifies closed manifolds carrying real Killing spinors in odd dimensions. We summarize the results in the table below.

	\begin{table}[H]
		\centering
	\begin{tabular}{|c|c|l|}
		\hline
		\textbf{Dimension of $M$} & \textbf{Type} & \textbf{Manifolds} \\
		\hline
		\rule{0pt}{12pt}
		$n$ arbitrary & $(2^{[\frac{n}{2}]},2^{[\frac{n}{2}]})$ & $\mathbb{S}^n$\\
		\hline
		\rule{0pt}{12pt}
		$n=2m-1$, $m\geq3$ & $(1,0)$ & Sasaki-Einstein manifold\\
		\hline
		\rule{0pt}{12pt}
		$n=4m-1$, $m\geq3$ & $(2,0)$  &  Sasaki-Einstein \\
		 &  &   without a Sasaki $3$-structure \\
		 & $(m+1,0)$& $M$ carries a Sasaki $3$-structure  \\ 
		\hline
		\rule{0pt}{12pt}
		$n=6$ & $(1,1)$ & $M$ nearly Kähler of constant type $1$ \\
		\hline
		\rule{0pt}{12pt}
		$n=7$ & $(1,0)$ & $M$ carries $3$-form $\phi$ such that $\nabla \phi=*\phi$ \\
		 & $(2,0)$ & $M$ carries Sasaki structure, \\
		  &  &  but no Sasaki $3$-structure \\
		 & $(3,0)$ & $M$ carries Sasaki $3$-structure \\
		\hline
		
	\end{tabular}
		\caption{Classification of closed manifolds carrying real Killing spinors based on the results in \cite{baer}.}\label{tab:characterization}
	\end{table}
We will omit details on the kind of more exceptional manifolds in Table \ref{tab:characterization} and refer the interested reader to the book \cite{boyer}. For our purpose, the most important consequence in the table above is that in dimension $n=3$, we directly obtain $(M,g)\cong (\mathbb{S}^3,g_{\mathbb{S}^3})$. Thus, in three dimensions equality in \eqref{eq:inequYamspin} holds if and  only if $(M,g)$ is conformally equivalent to $(\mathbb{S}^3,g_{\mathbb{S}^3})$. \\
Finally, we note that for even $n\neq 6$, we can not attain equality in \eqref{eq:inequYamspin} on $\mathbb{S}^n$, as by \cite[Theorem 1]{FrankLoss2} equality on $\mathbb{R}^n$ is attained if and only if $n$ is odd. As already explained above, we know that any manifold on which equality in \eqref{eq:inequYamspin} is attained is conformally equivalent to a manifold carrying a nontrivial Killing spinor. Thus, by \cite[Theorem 1]{baer} we obtain that $(M,g)$ is already conformally equivalent to $(\mathbb{S}^n,g_{\mathbb{S}^n})$ for all $n\neq 6$ even. Thus, in any even dimension apart from $n=6$ we conclude that equality in \eqref{eq:inequYamspin} never holds.
}
	\bibliographystyle{plain}
	\bibliography{zeromodesrefs}

\begin{thebibliography}{10}

\bibitem{amm2}
Bernd Ammann.
\newblock A variational problem in conformal spin geometry.
\newblock Habilitationsschrift, 2003.

\bibitem{amm1}
Bernd Ammann and Nadine Gro{\ss}e.
\newblock Relations between threshold constants for {Y}amabe type bordism
  invariants.
\newblock {\em J. Geom. Anal.}, 26(4):2842--2882, 2016.

\bibitem{aub}
Thierry Aubin.
\newblock \'equations diff\'erentielles non lin\'eaires et probl\`eme de
  {Y}amabe concernant la courbure scalaire.
\newblock {\em J. Math. Pures Appl. (9)}, 55(3):269--296, 1976.

\bibitem{aub2}
Thierry Aubin.
\newblock Probl\`emes isop\'erim\'etriques et espaces de {S}obolev.
\newblock {\em J. Differential Geometry}, 11(4):573--598, 1976.

\bibitem{aubbook}
Thierry Aubin.
\newblock {\em Some nonlinear problems in {R}iemannian geometry}.
\newblock Springer Monographs in Mathematics. Springer-Verlag, Berlin, 1998.

\bibitem{friedrich_baum}
Helga Baum and Thomas Friedrich.
\newblock Eigenvalues of the {D}irac operator, twistors and {K}illing spinors
  on {R}iemannian manifolds.
\newblock In {\em Clifford algebras and spinor structures}, volume 321 of {\em
  Math. Appl.}, pages 243--256. Kluwer Acad. Publ., Dordrecht, 1995.

\bibitem{boyer}
Charles~P. Boyer and Krzysztof Galicki.
\newblock {\em Sasakian geometry}.
\newblock Oxford Mathematical Monographs. Oxford University Press, Oxford,
  2008.

\bibitem{baer}
Christian Bär.
\newblock Real {K}illing spinors and holonomy.
\newblock {\em Comm. Math. Phys.}, 154(3):509--521, 1993.

\bibitem{FrankLoss}
Rupert~L. Frank and Michael Loss.
\newblock Which magnetic fields support a zero mode?
\newblock {\em J. Reine Angew. Math.}, 788:1--36, 2022.

\bibitem{FrankLoss2}
Rupert~L. Frank and Michael Loss.
\newblock A sharp criterion for zero modes of the dirac equation.
\newblock {\em Journal of the European Mathematical Society, published online
  first}, June 2024.
\newblock Available at \url{https://doi.org/10.4171/jems/1475}.

\bibitem{friedrich}
Thomas Friedrich.
\newblock {\em Dirac-{O}peratoren in der {R}iemannschen {G}eometrie}.
\newblock Advanced Lectures in Mathematics. Friedr. Vieweg \& Sohn,
  Braunschweig, 1997.
\newblock Mit einem Ausblick auf die Seiberg-Witten-Theorie. [With an outlook
  on Seiberg-Witten theory].

\bibitem{gro2}
Nadine Gro{\ss}e.
\newblock Solutions of the equation of a spinorial {Y}amabe-type problem on
  manifolds of bounded geometry.
\newblock {\em Comm. Partial Differential Equations}, 37(1):58--76, 2012.

\bibitem{gro}
Nadine Gro{\ss}e.
\newblock The {Y}amabe equation on manifolds of bounded geometry.
\newblock {\em Comm. Anal. Geom.}, 21(5):957--978, 2013.

\bibitem{gro1}
Nadine Gro{\ss}e and Cornelia Schneider.
\newblock Sobolev spaces on {R}iemannian manifolds with bounded geometry:
  general coordinates and traces.
\newblock {\em Math. Nachr.}, 286(16):1586--1613, 2013.

\bibitem{hebey}
Emmanuel Hebey.
\newblock {\em Nonlinear analysis on manifolds: {S}obolev spaces and
  inequalities}, volume~5 of {\em Courant Lecture Notes in Mathematics}.
\newblock New York University, Courant Institute of Mathematical Sciences, New
  York; American Mathematical Society, Providence, RI, 1999.

\bibitem{hij}
Oussama Hijazi.
\newblock A conformal lower bound for the smallest eigenvalue of the {D}irac
  operator and {K}illing spinors.
\newblock {\em Comm. Math. Phys.}, 104(1):151--162, 1986.

\bibitem{jin}
Zhi~Ren Jin.
\newblock A counterexample to the {Y}amabe problem for complete noncompact
  manifolds.
\newblock In {\em Partial differential equations ({T}ianjin, 1986)}, volume
  1306 of {\em Lecture Notes in Math.}, pages 93--101. Springer, Berlin, 1988.

\bibitem{lepa}
John~M. Lee and Thomas~H. Parker.
\newblock The {Y}amabe problem.
\newblock {\em Bull. Amer. Math. Soc. (N.S.)}, 17(1):37--91, 1987.

\bibitem{lossyau}
Michael Loss and Horng-Tzer Yau.
\newblock Stability of {C}oulomb systems with magnetic fields. {III}. {Z}ero
  energy bound states of the {P}auli operator.
\newblock {\em Comm. Math. Phys.}, 104(2):283--290, 1986.

\bibitem{roe}
John Roe.
\newblock {\em Elliptic operators, topology and asymptotic methods}, volume 179
  of {\em Pitman Research Notes in Mathematics Series}.
\newblock Longman Scientific \& Technical, Harlow; copublished in the United
  States with John Wiley \& Sons, Inc., New York, 1988.

\bibitem{schoe}
Richard Schoen.
\newblock Conformal deformation of a {R}iemannian metric to constant scalar
  curvature.
\newblock {\em J. Differential Geom.}, 20(2):479--495, 1984.

\bibitem{shu}
M.~A. Shubin.
\newblock Spectral theory of elliptic operators on non-compact manifolds.
\newblock In {\em M\'ethodes semi-classiques. Vol. 1. \'Ecole d'\'et\'e
  (Nantes, juin 1991)}, pages 35--108. Paris: Soci{\'e}t{\'e} Math{\'e}matique
  de France, 1992.

\bibitem{stru}
Michael Struwe.
\newblock {\em Variational methods}, volume~34 of {\em Ergebnisse der
  Mathematik und ihrer Grenzgebiete. 3. Folge. A Series of Modern Surveys in
  Mathematics [Results in Mathematics and Related Areas. 3rd Series. A Series
  of Modern Surveys in Mathematics]}.
\newblock Springer-Verlag, Berlin, fourth edition, 2008.
\newblock Applications to nonlinear partial differential equations and
  Hamiltonian systems.

\bibitem{tal}
Giorgio Talenti.
\newblock Best constant in {S}obolev inequality.
\newblock {\em Ann. Mat. Pura Appl. (4)}, 110:353--372, 1976.

\bibitem{tru}
Neil~S. Trudinger.
\newblock Remarks concerning the conformal deformation of {R}iemannian
  structures on compact manifolds.
\newblock {\em Ann. Scuola Norm. Sup. Pisa Cl. Sci. (3)}, 22:265--274, 1968.

\bibitem{yam}
Hidehiko Yamabe.
\newblock On a deformation of {R}iemannian structures on compact manifolds.
\newblock {\em Osaka Math. J.}, 12:21--37, 1960.

\end{thebibliography}
\end{document}